\documentclass{article}
\usepackage{ifthen}
\usepackage{graphicx}

\begin{document}

\title{Proof of the cases $p \leq 7$ of
the Lieb-Seiringer formulation of the Bessis-Moussa-Villani
conjecture}

\author{Daniel~H\"{a}gele\footnote{Spektroskopie der kondensierten Materie,
Ruhr--Universit\"at--Bochum,
              Universit\"atsstra\ss e 150,
              44780 Bochum, Germany, daniel.haegele@ruhr--uni--bochum.de}}
              \maketitle

\begin{abstract}
It is shown that the polynomial $\lambda(t) = {\rm Tr}[(A +
tB)^p]$ has nonnegative coefficients when $p \leq 7$ and A and B
are any two complex positive semidefinite $n \times n$ matrices
with arbitrary $n$. This proofs a general nontrivial case of the
Lieb-Seiringer formulation of the Bessis-Moussa-Villani conjecture
which is a long standing problem in theoretical physics.
\end{abstract}

In 1975 Bessis, Moussa, and Villani (BMV) stated the following,
now widely known,  conjecture. For arbitrary hermitian matrices
$G$ and $H$ the function $t \rightarrow {\rm tr} \exp(G + i t H)$
is the Fourier transform of a positive measure
\cite{bessisJMP75,moussaRMP00}. The conjecture is highly relevant
in the context of quantum mechanical partition functions and their
derivatives. However, as of today, no complete proof of the
conjecture is known, despite "a lot of work, some by prominent
mathematical physicists" as Reinhard F. Werner puts it on his
homepage \cite{wernerWEB96}.

Recently, Lieb and Seiringer found an equivalent formulation of
the BMV conjecture \cite{liebJSP04} which triggered new attempts
to prove it \cite{hillarJSP05,hillarMATH07}:

 {\bf
Conjecture} (Bessis--Moussa--Villani) The polynomial $\lambda(t) =
{\rm Tr} [(A+ t B)^p]$ has nonnegative coefficients whenever $A$
and $B$ are $n$-by-$n$ positive semidefinite matrices.

So far, only the trivial cases for $p \leq 5$ (for references see
\cite{hillarMATH07}) along with a single non-trivial case $p = 6$,
and matrix dimension $n = 3$ have been proven \cite{hillarJSP05}.

In the following I give an explicit elementary proof for the case
$p=7$ with no restrictions on the size of the matrices $A$ and
$B$. This means that for the first time a general non-trivial
result on the BMV conjecture for matrix dimensions larger than
three could be obtained. According to a theorem of C. J. Hillar
(corollary 1.8 in \cite{hillarMATH07}), the proof of $p=7$ implies
that all cases with $p \leq 7$ are proven. Especially the case $p
= 6$ and $n = 3$ is covered, which previously had been proven in a
completely different way by intricate computer algebra
\cite{hillarJSP05}.

The coefficient of $t^r$ in $\lambda(t)$ is the trace of
$S_{p,r}(A,B)$, the sum of all words of length $p$ in $A$ and $B$,
in which $r$ $B$'s appear (sometimes called the $r$-th Hurwitz
product of $A$ and $B$).

{\bf Theorem} The polynomial $\lambda(t) = {\rm Tr} [(A+ t B)^p]$
has non-negative coefficients when $p = 7$ and $A$ and $B$ are any
two $n$--by--$n$ positive semidefinite matrices.

{\it Proof.} Consider the third coefficient ${\rm Tr}S_{7,3}$. The
sum $S_{7,3}$ consists of 35 products where $A$ appears four times
and $B$ three times in all possible permutations.
Because of the cyclicity of the trace, we can write ${\rm
Tr}S_{7,3}$ as a sum of only five traces
\begin{equation}
  {\rm Tr}S_{7,3} = 7 (T_1 +  T_2 +  T_3 +  T_4 + T_5), \label{Tsum7}
\end{equation}
where
\begin{eqnarray}
  T_1 &=& {\rm Tr}(AAAABBB)  \nonumber\\
  T_2 &= &{\rm Tr}(AAABABB) \nonumber \\
  T_3 &=& {\rm Tr}(AAABBAB)  \nonumber\\
  T_4& =& {\rm Tr}(AABAABB) \nonumber \\
  T_5 & = & {\rm Tr}(AABABAB). \nonumber
\end{eqnarray}

To prove that ${\rm Tr}S_{7,3}$ is nonnegative I will show that it
can be written as a sum of traces of the form $Q = {\rm Tr} C^*
C$, where $C$ is a suitable chosen complex matrices. The relation
${\rm Tr} C^* C \geq 0$ holds for any $C$ because any matrix of
the form $C^* C$ is  positive semidefinite. Consider
\begin{eqnarray} C & = & c_1 bAAB + c_2 bABA + c_3 bBAA, \nonumber
\end{eqnarray}
where the $c_j$s are complex coefficients and $a$ and $b$ are
hermitian $n$--by--$n$ matrices with the properties $a^2 = A$ and
$b^2 = B$. The matrices $a$ and $b$ exist because $A$ and $B$ are
positive semidefinite and therefore also hermitian. We obtain
\begin{eqnarray}
  Q(c_1,c_2,c_3) & = & {\rm Tr} CC^* \nonumber\\
   & = & {\rm Tr}(c_1 c_1^* bAABBAAb + c_1 c_2^* bAABABAb + c_1 c_3^* bAABAABb \nonumber \\
  & & +  c_2 c_1^* bABABAAb + c_2 c_2^* bABAABAb + c_2 c_3^*
  bABAAABb
  \nonumber \\
  & & + c_3 c_1^* bBAABAAb + c_3 c_2^* bBAAABAb + c_3 c_3^* bBAAAABb)
  \nonumber \\
  & = &
c_1 c_1 T_4 + c_1 c_2^* T_5 + c_1 c_3^*
  T_4 \nonumber \\
& & + c_2 c_1^* T_5 + c_2 c_2^* T_5 + c_2 c_3^* T_3 \nonumber \\
& &  + c_3 c_1^* T_4 + c_3 c_2^* T_2 + c_3 c_3^* T_1  \nonumber
\geq 0.
\end{eqnarray}
By comparison with (\ref{Tsum7}) we find
\begin{eqnarray}
 {\rm Tr}
 S_{7,3} & =  & 7 Q(1,0,0) + 7 Q(0,1,1) \geq 0 ,
   \nonumber
\end{eqnarray}
which concludes the proof of the third coefficient being
nonnegative.\\
Similarly to the proof above, the nonnegativity of the $0$th,
$1$st and $2$nd coefficient can be proven:
\begin{eqnarray}
 {\rm Tr}S_{7,2} & = &  7 {\rm Tr}(BAAAAAB + ABAAABA + AABABAA) \nonumber
\\
 & = & 7 {\rm Tr} (BAAa * H.C. + ABAa * H.C. + AABa * H.C.) \geq
 0, \nonumber \\
  {\rm Tr}S_{7,1} & = & 7 {\rm Tr} (AAABAAA) \nonumber \\
          & = & 7 {\rm Tr} (AAAb * H.C.) \geq 0, \nonumber \\
  {\rm Tr}S_{7,0} & = &  {\rm Tr} (AAAAAAA) \nonumber \\
          & = &  {\rm Tr} (AAAa * H.C.) \geq 0, \nonumber
\end{eqnarray}
where H.C. denotes the hermitian conjugate of the first factor in
the product. The cases ${\rm Tr} S_{7,4}$,${\rm Tr} S_{7,5}$,${\rm
Tr} S_{7,6}$ and $S_{7,7}$ are, after exchanging $A$ and $B$,
identical with the cases ${\rm Tr} S_{7,3}$,${\rm Tr} S_{7,2}$,
${\rm Tr} S_{7,1}$, and $S_{7,0}$, respectively and are therefore
covered by the proof above. This completes the proof of the
theorem.

{\bf Corollary} The polynomial $\lambda(t) = {\rm Tr} [(A+ t
B)^p]$ has non-negative coefficients when $p \leq 7$ and $A$ and
$B$ are any two $n$--by--$n$ positive semidefinite matrices.

 {\it Proof.} The corollary follows immediately from a theorem of C. J.
Hillar:

 {\bf Theorem} (corollary 1.8 in \cite{hillarMATH07}). If the
Bessis-Moussa-Villani conjecture (see above) is true for some
$p_0$, then it is also true
for all $p < p_0$.\\

{\it The case $p = 6$.} One may ask if it is possible to find a
direct way of proving the case $p=6$ in a fashion similar to the
proof shown above. This is not possible as I will show below.
Analogously to the case $p = 7$, we find after some algebra
\begin{eqnarray}
  {\rm Tr}S_{6,3} & =&  6 T_1 + 6 T_2 + 6 T_3 + 2 T_4. \label{Tsum}
\end{eqnarray}
where
\begin{eqnarray}
  T_1 &=& {\rm Tr}(AAABBB)  \nonumber\\
  T_2 &= &{\rm Tr}(AABABB)  \nonumber \\
  T_3 &=& {\rm Tr}(AABBAB)  \nonumber\\
  T_4& =& {\rm Tr}(ABABAB). \nonumber
\end{eqnarray}
Using
\begin{eqnarray} C & = &c_1 \, aABb + c_2 \,aBAb  \nonumber
\end{eqnarray}
we obtain
\begin{eqnarray}
  Q(c_1, c_2) & = &  {\rm Tr} C C^* \nonumber \\
  &  = &  {\rm Tr} (c_1 c_1^* aABbbBAa + c_1 c_2^* aABbbABa
  \nonumber
  \\
   & & c_2 c_1^* aBAbbBAa + c_2 c_2^* aBAbbABa ) \geq 0.
   \nonumber
\end{eqnarray}
and
\begin{eqnarray}
  Q(c_1,c_2) & = & c_1 c_1^* T_1 + c_1 c_2^* T_3 +
  c_2 c_1^* T_2 + c_2 c_2^* T_4 \geq 0.  \nonumber
\end{eqnarray}
Here it is impossible to find a way to write ${\rm Tr} S_{6,3}$ as
a sum of $Q(c_1,c_2)$: Suppose there is a solution, then from
\begin{eqnarray}
  {\rm Tr} S_{6,3} & = & \sum_l Q(c_{1,l},c_{2,l}) \nonumber \\
    & = & \sum_l c_{1,l}c_{1,l}^* T_1 + \sum_l c_{1,l}c_{2,l}^* T_3 +  \sum_l c_{2,l}c_{1,l}^*
    T_2 + \sum_l c_{2,l}c_{2,l}^* T_4  \nonumber
\end{eqnarray}
we find by comparison with eqn. (\ref{Tsum}) $\sum_j
c_{1,l}c_{1,l} = 6$, $\sum_l c_{1,l}c_{2,l}^* = 6$, $\sum_l
c_{2,l}c_{1,l}^* = 6$, and $\sum_l c_{2,l}c_{2,l}^* = 2$. Since
generally $(c_{1,l} - c_{2,l})(c_{1,l} - c_{2,l})^* \geq 0$ we
find
\begin{eqnarray}
 \sum_l
c_{1,l}c_{1,l}^* + \sum_l c_{2,l}c_{2,l}^* & \geq & \sum_l
c_{1,l}c_{2,l}^*
     + \sum_l c_{2,l}c_{1,l}^* \nonumber \\
     6 + 3 & \geq  & 6 +6,  \nonumber
\end{eqnarray}
which contradicts our assumption that ${\rm Tr} S_{6,3}$ can be
written as sum of $Q$s. Note that the $C$ was chosen to yield the
right number of $A$s and $B$s in the sum. We speculate that a more
general ansatz that relaxes this restriction for $C$, i.e. $C_1
=c_1 \, aABb + c_2 \,aBAb + c_3\, aAAb + c_4\, aBBb$ and $C_2 =
d_1 aAAa + d_2 aABa + ...$, $C_3 = ...$, ... may lead to a direct
proof.

In conclusion, I gave an elementary proof of the nontrivial case
$p=7$ of the Lieb-Seiringer formulation of the
Bessis-Moussa-Villani conjecture for $n \times n$-matrices.
Although a proof of $p=7$ implies the $p = 6$ case, the way of
proving the $p=7$ case did not apply analogously. Attempts to
prove cases with $p > 7$ by methods similar to those presented
here will be more complex, but should nevertheless be undertaken
in the future using computer algebra.

\end{document}